\documentclass[11pt,a4paper]{article}

\usepackage[top=3cm, left=3.5cm, right=3cm, bottom=3cm]{geometry}
\usepackage[activate]{pdfcprot}

\usepackage{hyperref}

\hypersetup{%
 pdftitle = {Regularization of statistical inverse problems and the Bakushinskii veto}
 pdfauthor = {Saskia Becker},
 pdfkeywords = {statistics, inverse problem, stochastical noise, regularization, Bakushinskii},
 pdfstartview = {FitH},
 colorlinks = {false},
 breaklinks = {true},
 pdfborder = {0 0 0},
 linkcolor = {blue},
  plainpages = {false},
}

\usepackage{paralist}

\usepackage{pslatex}

\makeatletter 
\renewcommand{\section}{\@startsection
{section}
{1}
{0mm}
{2\baselineskip}
{\baselineskip}
{\normalfont\large\scshape}}
\makeatother

\makeatletter 
\renewcommand{\subsection}{\@startsection
{subsection}
{2}
{0mm}
{2\baselineskip}
{1\baselineskip}
{\normalfont\normalsize\itshape}}
\makeatother

\usepackage[ansinew]{inputenc}

\usepackage{fancyhdr}

\usepackage[square,numbers]{natbib}

\usepackage[nottoc]{tocbibind} 

\usepackage{amsmath}
\usepackage{amsfonts}
\usepackage{graphicx}
\graphicspath{{Graphiken/}}
\usepackage{subfig}

\usepackage{amsthm}
\numberwithin{equation}{section}
\theoremstyle{definition}
\newtheorem{Def}{Definition}[section]
\newtheorem{Not}[Def]{Notation}
\newtheorem{Ass}{Assumption}
\theoremstyle{plain}
\newtheorem{Thm}[Def]{Theorem}
\newtheorem{Prop}[Def]{Proposition}
\newtheorem{Lem}[Def]{Lemma}
\theoremstyle{remark}
\newtheorem{Conc}[Def]{Conclusion}
\newtheorem{Rem}[Def]{Remark}
\newtheorem{Ex}[Def]{Example}

\begin{document}

\title{\textbf{Regularization of statistical inverse problems\\ and the Bakushinski{\u\i} veto}}
\author{Saskia Becker\\
  \small Weierstra\ss{} Institute for Applied Analysis and Stochastics\\
  \small Mohrenstra\ss{}e 39, 10117 Berlin, Germany\\
  \small saskia.becker@wias-berlin.de}

\thispagestyle{empty}
\maketitle

\begin{quote}
\textbf{Abstract.} In the deterministic context Bakushinski{\u\i}'s theorem excludes the existence of purely data driven convergent regularization for ill-posed problems. We will prove in the present work that in the statistical setting we can either construct a counter example or develop an equivalent formulation depending on the considered class of probability distributions. Hence, Bakushinski{\u\i}'s theorem does not generalize to the statistical context, although this has often been assumed in the past. To arrive at this conclusion, we will deduce from the classic theory new concepts for a general study of statistical inverse problems and perform a systematic clarification of the key ideas of statistical regularization. 
\end{quote}


\section{Introduction}

%

We consider statistical inverse problems, where an unknown signal $x$ should be reconstructed from indirect noisy measurements $ y_{\mathrm{noise}} = Tx + \mathrm{noise}$. The problem is assumed to be ill-posed, i.e. the operator $T$ is not continuously invertible such that we can only approximate the signal. In classic inverse problems the noise is supposed to be deterministic and bounded. Nevertheless it is well-known that various applications cannot be modeled appropriately in this way. Therefore, stochastic models have been introduced, where the noise is taken as random variable or stochastical process \cite{MR2505875, MR2361904, hofinger, MR2503292}. In some studies, e.g. \cite{MR2158113, hofinger, MR859375}, not only the noise but also the operator or the signal are stochastic. 

\vspace{11 pt}
In both the deterministic and the stochastic setting one crucial point is the knowledge of the noise level which is often not available in application. However, the Bakushinski{\u\i} veto \cite{MR760252} states for classic inverse problems the equivalence of the ill-posedness of the problem and the nonexistence of purely data driven reconstruction methods, for which the approximated solution tends to the exact signal $x$ when the noise vanishes. This theorem is of particular importance since it constitutes the need of supplemental information, as for instance the noise level.

\vspace{11 pt}
For statistical inverse problems the situation is ambiguous as we will discuss in the paper at hand. To study the existence of such reconstruction methods we need explicit definitions of the involved objects. While an extensive theory for classic inverse problems has been developed \cite{MR1408680, MR742928, MR859375}, only selected aspects of statistical inverse problems have been analyzed so far. Additional difficulties, arising from the possible unboundedness of stochastic noise, are the need of new error and convergence criteria \cite{MR2361904, MR1929274, MR2503292, MR874480}. Cavalier explained in \cite{MR2421941} how concepts of nonparametric statistics, e.g. the white noise model, risk estimation and model selection, can be applied to inverse problems. We will proceed in reverse by studying how the key ideas of the classic inversion theory have to be modified for beeing suitable for a statistical setting.

\vspace{11 pt}
First of all we give a brief recapitulation of the classic regularization theory, in which we suggest in particular a reduction of the usually required convergence properties. Our statistical setting is introduced in section \ref{sec:setting} being followed by the presentation of the main concepts and central definition in section \ref{sec:kindsConv}. There we propose to link the noise to the asymptotic of the noise level, which will turn out to be the deciding idea for definition \ref{def:stKonvRV} of convergent statistical regularization methods and our main result stated in section \ref{sec:stBak}: We prove an equivalent formulation and give a counterexample to Bakushinski{\u\i}'s theorem depending on the considered class of probability distributions.


\section{Classic inverse problems}\label{sec:clIP}

We consider the usual setting of classic inverse problems. Let $\mathbb{H}_{1}$ and $\mathbb{H}_{2}$ denote separable Hilbert spaces with scalar products $\langle .,. \rangle_{\mathbb{H}_{i}}$ and the induced norms $\Vert . \Vert_{\mathbb{H}_{i}}$, $i=1,2$. Further let $T: \mathbb{H}_{1} \rightarrow \mathbb{H}_{2}$ be a linear, compact and bounded operator with a nonclosed range $\mathcal{R}(T)$. We are interested in the problem
\begin{equation}\label{eq:detyd}
y_{\delta} = Tx + \delta \xi,
\end{equation}
where $x \in \mathbb{H}_{1}$ denotes the unknown signal, $\delta > 0$ is the noise level and the normalized noise $\xi \in \mathbb{H}_{2}$ satisfies $\Vert \xi \Vert_{\mathbb{H}_{2}} \leq 1$. With $\ker(T)^{\perp}$ as orthogonal complement of the kernel of $T$ we can define the generalized inverse $T^{+}$ as the linear extension of the inverse of $ T\vert_{\ker(T)^{\perp}} $. A motivation and some properties of the generalized inverse can be found e.g. in \cite{MR1408680}. Since the range of $T$ is assumed to be nonclosed, $T^{+}$ is discontinuous and $x^{+} := T^{+} y  \in \mathbb{H}_{1}$ has to be regularized.

\vspace{11 pt}
In the following subsection we will not present the common definition of (convergent) regularization methods given in \cite{MR1408680}, but the definitions introduced by Hofmann and Math\'{e} in \cite{MR2318806}. Research has shown that purely data driven regularization methods can yield remarkably good results, see for instance \cite{MR2175028, MR1408680}, although these methods are not convergent as the Bakushinski{\u\i} veto proves. This teaches us to distinguish convergent and arbitrary regularization schemes as is done in the following approach.


\subsection{Linear and convergent regularization schemes}\label{sec:LinStReg}

\begin{Not}[Singular value decomposition (SVD) of $T$ \cite{MR1408680}] 
Let $ \left\{ \left( s_{j}; v_{j}, u_{j} \right) \right\}_{j \in \mathbb{N}}$ denote the singular system of the operator $T$, where $\left\{ s_j \right\}_{j \in \mathbb{N}}$ is arranged in decreasing order with $ \underset{j \rightarrow \infty}{\lim} s_{j} = 0$. The following series expansion holds:
\[
 Tx = \sum_{j \geq 1} s_{j} \left\langle x, v_{j} \right\rangle_{\mathbb{H}_{1}} u_{j}
\]
\end{Not}

\begin{Def}[Linear regularization \cite{MR2318806}]\label{def:klRV}
A family $ F := \left\{ F_{\alpha} \right\}_{\alpha>0} $ of linear and bounded operators $F_{\alpha}: \left[ 0, \left\| T \right\|^{2} \right] \rightarrow \mathbb{R}$ is called regularization (filter) if the following properties hold:
	\begin{compactenum}
	\item The associated bias family $\left\{ b_{\alpha} \right\}_{\alpha>0}$, where $b_{\alpha}(\vartheta) := 1 - \vartheta F_{\alpha}(\vartheta)$, converges pointwise to zero: $\underset{\alpha \rightarrow 0}{\lim} \, b_{\alpha}(s_{j}^2 ) = 0$ for all $s_j>0$.
	\item The bias family is uniformly bounded by some $\gamma_{0} > 0$, i.e. $\underset{\alpha \leq \alpha_{0}}{\sup} \;\underset{s_{j}>0}{\sup} \left| b_{\alpha}(s_j^2) \right| \leq \gamma_{0}$.
	\item There is a constant $\gamma_{*} > 0$ such that the parameter family can be normalized for all $\alpha \in (0,\infty)$ and $s_{j}>0$ by $ s_{j} \left| F_{\alpha}(s_{j}^2) \right| < \gamma_{*} / \sqrt{\alpha}$.
	\end{compactenum}
In this case, the family $ R := \left\{ R_{\alpha} \right\}_{\alpha > 0} $ of linear and bounded operators $ R_{\alpha} : \mathbb{H}_{2} \rightarrow \mathbb{H}_{1}$ with
	\begin{equation}\label{eq:klROperator} 
		R_{\alpha} y := x_{\alpha} := F_{\alpha}(T^{*}T)T^{*}y = \sum_{s_{j} > 0} F_{\alpha}(s_{j}^{2}) s_j \left\langle y, u_{j} \right\rangle_{\mathbb{H}_{1}} v_{j}, \, y \in \mathbb{H}_{2},
	\end{equation} 
is called linear regularization scheme (in short: regularization), where the last equation follows from the functional calculus described in \cite{MR1408680}.
\end{Def}

\begin{Not}
Below, we will use without further comments the notations
\[
F := \left\{ F_{\alpha} \right\}_{\alpha>0} \qquad \text{ and } \qquad R := \left\{ R_{\alpha} \right\}_{\alpha > 0}.
\]
\end{Not}

\begin{Ex}\label{ex:RVen}\hspace{1 pt}
The given definition is satisfied by many of the known linear regularization in terms of \cite{MR1408680} such as spectral cut-off, which is defined by
 	\[
 		F_{\alpha}(\vartheta) := \vartheta^{-1} \, \chi_{\left( \alpha, \left\| T \right\|^{2} \right) } (\vartheta)
 		\quad
 	\text{ such that }
 		x_{\alpha} = R_{\alpha}y = \sum_{s_{j}^{2} > \alpha} s_{j}^{-1} \left\langle y, u_{j} \right\rangle v_{j},
 	\]
 where $\chi$ denotes the indicator function, $\alpha, \vartheta \in \left( 0, \left\| T \right\|^{2} \right] $, and Tikhonov regularization with
 	\[
 		F_{\alpha}(\vartheta) := 1/(\alpha + \vartheta),
 	\text{ such that }
 		x_{\alpha} = R_{\alpha}y = (\alpha I + T^{*}T)^{-1} T^{*} y.
 	\]
 \end{Ex}

\begin{Rem}\label{rem:linRVen}\hspace{1 pt}
\begin{compactenum}
\item Later on we will require a stricter bound instead of property (3) of definition \ref{def:klRV}:
\begin{equation}\label{eq:Vainikko3}
	\sup_{0 < \vartheta \leq \left\| T \right\|^{2}} \vert F_{\alpha} (\vartheta) \vert \leq \tfrac{\gamma}{\alpha}, \, \gamma > 0.
\end{equation}
It is easy to show, that the given examples satisfy this property too. In \cite{MR859375} it is shown, that (3) follows if (2) and (\ref{eq:Vainikko3}) hold.
\item As generalization we could also require that the index family of $F$ is an arbitrary subset of the real numbers with at least one accumulation point, say $h \in \mathbb{R}$. Then property (1) has to be reformulated in the following way: $\underset{\alpha \rightarrow h }{\lim} \, b_{\alpha}(s_j^2) = 0$ for all $ s_j>0 $. We cannot skip it completely because it yields the following important proposition.
\end{compactenum}
\end{Rem}

\begin{Prop}[Pointwise convergence to $T^{+}$ \cite{matheVL}]\label{thm:R-pktw-T+}
Let $ R $ denote a linear regularization and $\mathcal{D}(T^{+})$ the domain of the generalized inverse $T^{+}$ of $T$.
\begin{compactenum}
	\item If $y \in \mathcal{D}(T^{+})$, then
	$
		\underset{\alpha \leq \alpha_{0}}{\sup} \left\| x_{\alpha} \right\|_{\mathbb{H}_{1}} < \infty \text{ and }
		x_{\alpha} = R_{\alpha}y \rightarrow T^{+}y \text{ when } \alpha \rightarrow 0.
	$
	\item If $y \notin \mathcal{D}(T^{+})$, then $\underset{\alpha \rightarrow 0}{\lim} \left\| x_{\alpha} \right\|_{\mathbb{H}_{1}} = \infty $.
\end{compactenum}
	In particular we get for all $y \in \mathbb{H}_{2}$ that $ \underset{\alpha \rightarrow 0}{\lim}\,  T R_{\alpha}y = TT^{+}y = Qy$, where $Q: \mathbb{H}_{2} \rightarrow \overline{\mathcal{R}(T)}$ denotes the orthogonal projection onto $\overline{\mathcal{R}(T)}$.
\end{Prop}

\begin{Rem}
A similar result can be found in \cite[proposition 3.6]{MR1408680}.
\end{Rem}

Convergence in general and especially convergence rates are established quality criteria for the comparison of regularization schemes. Normally, one claims that the regularized solution $x_{\alpha}$ should converge uniformly to the exact one, if the error tends to zero:

\begin{Def}[Parameter choice \cite{matheVL}]\label{def:klPW}\label{def:klKonvRV}
Let $ R $ denote a linear regularization scheme and $\alpha: (0, \infty) \times \mathbb{H}_{2}\rightarrow (0, \infty)$ a function. If for all $y \in \mathcal{D}(T^{+})$ it holds
\[
	 \lim_{\delta \rightarrow 0} \left( \sup \left\{ \alpha(\delta, y_{\delta}) : y_{\delta} \in \mathbb{H}_{2}, \left\| y - y_{\delta} \right\|_{\mathbb{H}_{2}} \leq \delta \right\} \right) = 0,
\] 
then $\alpha$ is called (classic) parameter choice. In particular we will say:
\begin{compactenum}
\item $\alpha$ is purely data driven or heuristic if it depends only on the data, i.e. $\alpha = \alpha(y_{\delta})$.
\item $\alpha$ is (classic) convergent w.r.t. $R$ if for all $y \in \mathcal{D}(T^{+})$ it holds
\[
	\lim_{\delta \rightarrow 0} \left( \sup \left\{ \left\| T^{+}y - R_{\alpha (\delta, y_{\delta}) } y_{\delta} \right\|_{\mathbb{H}_{1}} : y_{\delta} \in \mathbb{H}_{2}, \left\| y - y_{\delta} \right\|_{\mathbb{H}_{2}} \leq \delta \right\} \right) = 0.
\]
\end{compactenum}
The pair $(R, \alpha )$ of a linear regularization $R $ and a parameter choice $\alpha$ is called (classic) convergent regularization method of $T^{+}$ if $\alpha$ is convergent w.r.t. $R$.
\end{Def}

\begin{Not}\hfill
\begin{compactitem}
\item Here, we applied the usual error criterion for classic inverse problems:
\[
e(R, \alpha, x, \delta) :=  \sup \left\{ \left\| T^{+}y - R_{\alpha (\delta, y_{\delta}) } y_{\delta} \right\|_{\mathbb{H}_{1}} : y_{\delta} \in \mathbb{H}_{2}, \left\| y - y_{\delta} \right\|_{\mathbb{H}_{2}} \leq \delta \right\} \text{, where } y = Tx.
\]
\item Many parameter choice strategies depend on the applied regularization scheme $R$ which is why we should write $\alpha(R, \delta, y_{\delta})$. However, we will use for simplicity $\alpha(\delta, y_{\delta})$ instead.
\end{compactitem}
\end{Not}

\begin{Ex}\label{ex:detPW} 
The discrepancy principle \cite{MR1408680, MR0208819} is a good example of a parameter choice which is very common for classic inverse problems but cannot be applied in the statistical setting as we will explain in remark \ref{rem:DPinStIP}. It chooses the regularization parameter for a given regularization scheme $R$ and a fixed constant $\tau>1$ by setting
\[
	\alpha_{*} := \sup \left\{ \alpha \leq \left\| A^{*} A \right\|, \left\| A R_{\alpha} y_{\delta} - y_{\delta} \right\| \leq \tau \delta \right\}.
\]
Therein and in most of the established convergent methods the knowledge of the noise level $\delta$ is needed. In contrast, the quasi-solution of Ivanov \cite{MR2010817} yields convergent regularization assuming instead of that an upper bound for the norm of the solution $x_\alpha$. Well-known purely data driven parameter choices are the L-curve criterion of Hansen \cite{MR1193012}, the generalized cross-validation of Wahba \cite{MR1045442} and quasi-optimality 
\cite{MR0455365}.
\end{Ex}

\begin{Thm}[Bakushinski{\u\i} veto \cite{MR760252, MR1408680}]\label{thm:detBak}
A purely data driven (classic) convergent regularization method exists if and only if the generalized inverse $T^{+}$ is continuous.
\end{Thm}

\begin{proof}[Proof sketch]
With a purely data driven (classic) convergent regularization method $(R, \alpha)$ we get necessarily for exact data that $T^{+} y = R_{\alpha(y)} y$ for all $y \in  D(T^{+})$ such that for arbitrary sequences $\left\{y_{n}\right\}_{n\in\mathbb{N}} \subseteq  D(T^{+})$ with $ \underset{n\rightarrow \infty}{\lim} y_{n} = y$ it holds  $ \underset{n\rightarrow \infty}{\lim} T^{+} y_{n} = \underset{n\rightarrow \infty}{\lim} R_{\alpha(y_{n})}y_{n} = T^{+} y$, which yields the well-posedness of the problem.
\end{proof}


\subsection{Reduction of the requirements}

In the statistical setting we cannot require uniform convergence as we do in the deterministic context since the noise may be unbounded. The resulting question is, if for classic inverse problems the convergence criterion could also be diminished. We want to ensure that the approximated solution of the problem converges to the exact one if the noise tends to zero. But for that purpose we do not need to include the supremum as is done in definition \ref{def:klPW}.  It is only a technical simplification. Additionally we want to skip the requirement that the function $\alpha$ has to converge to zero if the noise vanishes. In fact, it is unimportant how $\alpha$ behaves as long as (\ref{eq:gKonRven}) is satisfied.


\begin{Def}[Generally convergent regularization]\label{def:konvRV}
The pair $(R, \alpha )$ of a linear regularization $R$ and a function $\alpha: (0, \infty) \times \mathbb{H}_{2}\rightarrow (0, \infty)$ is called (generally) convergent regularization of $T^{+}$ if the regularized solution converges in the following sense to the exact one: For all $ \left\{ y^{(k)} \right\}_{k \geq 1} $ with $y^{(k)} := y + \delta^{(k)} \xi^{(k)} $, $\delta^{(k)}> 0$, $ \left\| \xi^{(k)} \right\|_{\mathbb{H}_{2}} \leq 1 $ and $\underset{k \rightarrow \infty}{\lim} \delta^{(k)} = 0$ we have
\begin{equation}\label{eq:gKonRven}
	 \lim_{k \rightarrow \infty} \left\| T^{+}y - R_{\alpha (\delta^{(k)}, y^{(k)}) } y^{(k)} \right\|_{\mathbb{H}_{1}} = 0.
\end{equation}
\end{Def}


\begin{Rem}\label{rem:pointw}
In order to achieve an easier notation, one could be tempted to claim only pointwise convergence. But this would mean to fix the noise and vary only the noise level, which forms a considerable and unrealistic restriction.
\end{Rem}

\begin{Conc}[The Bakushinski{\u\i} veto for general methods]\label{rem:detBak}
As the supremum is not necessary for the proof of theorem \ref{thm:detBak}, an equivalent formulation can be varified analoguosly for generally convergent regularization.
\end{Conc}


\section{Statistical inverse problems}\label{sec:stIP}

In this section we provide new concepts for a general study of statistical inverse problems. As main idea we link the noise to the asymptotic of the noise level varying its probability distribution.


\subsection{Statistical setting}\label{sec:setting}

In recent publications about statistical inverse problems one can find two models of stochastical noise, random variables \cite{hofinger, MR2503292} and Hilbert-space processes \cite{MR2505875, MR2361904}. As every Hilbert-space valued random variable with finite second moment can be identified with a Hilbert-space process, we will concentrate mostly on the latter. 

\begin{Def}
A Hilbert-space process is a linear and continuous operator
\[
	\Xi: \mathbb{H}_{2} \rightarrow L^{2}(\Omega, \mathcal{F}, \mathbb{P}), \; v \mapsto \Xi v =: \left\langle \Xi, v \right\rangle_{\mathbb{H}_{2}},
\]
where $(\Omega, \mathcal{F}, \mathbb{P})$ denotes a probability space, $\mathcal{B}_{\mathcal{T}}$ the Borel-$\sigma$-algebra generated by the topological space $\mathcal{T}$ and
\[
	L^{2}(\Omega, \mathcal{F}, \mathbb{P}) := \left\{ Z: (\Omega, \mathcal{F}, \mathbb{P} ) \rightarrow (\mathbb{R}, \mathcal{B}_{\mathbb{R}}) \text{ square-integrable random variable} \right\}.
\]
\end{Def}

\begin{Def}
The covariance $ \mathrm{Cov}_{\Xi}: \mathbb{H}_{2} \rightarrow \mathbb{H}_{2} $ of a Hilbert-space process $\Xi$ is implicitly defined by
\[
	\left\langle \mathrm{Cov}_{\Xi} y_{1}, y_{2} \right\rangle_{\mathbb{H}_{2}} = \mathrm{Cov} \left( \left\langle \Xi, y_{1} \right\rangle_{\mathbb{H}_{2}} , \left\langle \Xi, y_{2} \right\rangle_{\mathbb{H}_{2}} \right), \, y_{1}, y_{2} \in \mathbb{H}_{2}.
\]
\end{Def}

Hence it is a bounded and linear operator.

\begin{Ex}\label{ex:GWR}
A centered Hilbert-space process $\Xi$ with the unit matrix as covarianceis called white noise process. In this case $\Xi$ is Gaussian if the associated random variables are Gaussian, i.e. if $	\left\langle \Xi, v \right\rangle_{\mathbb{H}_{2}} \sim \mathcal{N} \left( 0, \left\| v \right\|^{2}_{\mathbb{H}_{2}} \right)$.
Inverse problems with Gaussian white noise have been studied e.g. in \cite{MR2503292, MR2438944, MR2240642}.
\end{Ex}

\begin{Ass}\label{ass:Xi}
We assume $\Xi: \mathbb{H}_{2} \rightarrow L^{2}(\Omega, \mathcal{F}, \mathbb{P})$ to be a centered Hilbert-space process with $ \mathbb{E} \left[ \left\langle \Xi, v \right\rangle_{\mathbb{H}_{2}} \right] = 0 $ for all $ v \in \mathbb{H}_{2} $ and $ \left\| \mathrm{Cov}_{\Xi} \right\| < \infty $.
\end{Ass}

\begin{Not}[Observation model]\label{not:ObsModel}
Let $\Xi$ be as in assumption \ref{ass:Xi}. We consider the following abstract observation model:
\begin{equation}\label{eq:ObsModel}
 Y_{\delta} = y + \delta \Xi, \, \text{ where } y \in \mathcal{D}(T^{+}) \text{ and } \delta > 0.
\end{equation}
\end{Not}

\begin{Conc}
The realizations of $\Xi$ and thus of $Y_{\delta}$ do not have to be in $\mathbb{H}_{2}$ because $\Xi$ is only a weak random element of $\mathbb{H}_{2}$. As a consequence several basic concepts have to be revised:
\end{Conc}

\begin{Not}\label{not:InterprPXi}
We want to generalize the notation $\mathbb{P}^{\Xi}$ of image measures from random variables to Hilbert-space processes. Let $ \Xi $ be a Hilbert-space process. Then we interpret $\mathbb{P}^{\Xi}$ as the probability measure which is well-defined by its finite-dimensional marginal distributions on the space $(\mathbb{R}^{ \mathbb{H}_{2} }, (\mathcal{B}_{\mathbb{R}})^{ \otimes \mathbb{H}_{2} })$, where $ \mathbb{R}^{ \mathbb{H}_{2} } $ denotes the space of all functions $f: \mathbb{H}_{2} \rightarrow \mathbb{R}$ and $ (\mathcal{B}_{\mathbb{R}})^{ \otimes \mathbb{H}_{2} } $ denominates the associated $\sigma$-algebra. The existence and uniqueness of $\mathbb{P}^{\Xi}$ is ensured by the Kolmogorov extension theorem \cite{MR1083357}.
\end{Not}

\begin{Def}[Noise level]\label{def:Rauschpegel}
The definitions of the noise level of classic and statistical inverse problems differ significantly. The noise level $\delta$ of an inverse problem is defined as scale factor of the noise $\xi$ or accordingly $\Xi$, such that 
\begin{compactitem}
\item $ \left\| \xi \right\|_{\mathbb{H}_{2}} \leq 1$ and therefore $\left\| y - y_{\delta} \right\|_{\mathbb{H}_{2}} \leq \delta$ for all $\delta > 0$ if $y_{\delta}$ as in (\ref{eq:detyd}),
\item $ \mathbb{E} \left[ \left\| \Xi \right\|_{\mathbb{H}_{2}}^{2} \right] \leq 1 $ and therefore $\mathbb{E} \left[ \left\| y - Y_{\delta} \right\|_{\mathbb{H}_{2}}^{2} \right] \leq \delta^{2} $ for all $\delta > 0$ if $Y_{\delta} \in L^2(\Omega, \mathbb{H}_2)$,
\item $\left\| \mathrm{Cov}_{\Xi} \right\|^{1/2} \leq 1$ if $Y_{\delta}$ as in notation \ref{not:ObsModel}.
\end{compactitem}
\end{Def}

\begin{Rem}[Discussion]\label{rem:DPinStIP}
From a statistical point of view, only the third case is of interest. For instance, the discrepancy principle desribed in example \ref{ex:detPW} cannot be applied to observations with white noise since the term $\left\| A R_{\alpha} Y^{\delta} (\omega) - Y_{\delta} (\omega) \right\|_{\mathbb{H}_{2}}$ could be infinite. For observations with noise modelled as random variables it yields convergent methods by contrast. So, the second case is very close to the deterministic setting as we will support by proposition \ref{prop:TransfRVrv}.
\end{Rem}

For the deterministic context we defined the regularization operators between the observed Hilbert-spaces. The following notation allows us to apply them also to Hilbert-space processes:

\begin{Not}\label{not:RtXi-XiHP}
We observe a Hilbert-space process $\Xi: \mathbb{H}_{2} \rightarrow L^{2}(\Omega, \mathcal{F}, \mathbb{P})$ and a linear and bounded operator $R : \mathbb{H}_{2} \rightarrow \mathbb{H}_{1}$. Then, we will interpret the composition $R \: \Xi$ as a Hilbert-space process on $\mathbb{H}_{1}$, i.e. as
$ R \: \Xi : \mathbb{H}_{1} \rightarrow L^{2}(\Omega,  \mathcal{F}, \mathbb{P}) $ with $
v \mapsto R \: \Xi \: v =: \left\langle R \: \Xi, v \right\rangle_{\mathbb{H}_{1}} = \left\langle \Xi, R^{*} v \right\rangle_{\mathbb{H}_{2}}$.
\end{Not}

\begin{Rem}\label{rem:BemZuInterprRXi}	
$R \: \Xi$ is well-defined, since $\left\langle \Xi, R^{*} v \right\rangle_{\mathbb{H}_{2}} = \Xi (R^{*} v) \in L^{2}(\Omega, \mathcal{F}, \mathbb{P})$. The linearity of $R$ yields further that $R \: Y_{\delta} = R \: y + \delta R \: \Xi$.
\end{Rem}

As parameter choices do not have to be linear, we cannot interprete the term $\alpha \left( \delta, Y_{\delta} \right)$ in a similar way. That is why we will use, where necessary, the sequence space model, which was discussed for instance in \cite{MR2438944, MR2421941, MR2013911}:

\begin{Not}[Sequence space model]\label{not:SSModel}
Let $ \left\{ \left( s_{j}; v_{j}, u_{j} \right) \right\}$ denote the singular system of the operator $T$. The sequence space model is defined by
\begin{equation}\label{eq:Yjdelta}
	Y_{\delta} (\omega) := \left\{ Y_{\delta, j} (\omega) \right\}_{s_{j} > 0} \text{ with } Y_{\delta, j} (\omega)
	= \left\langle y, u_{j} \right\rangle_{\mathbb{H}_{2}} + \delta \left\langle \Xi, u_{j} \right\rangle_{\mathbb{H}_{2}} (\omega), \, \omega \in \Omega.
\end{equation}
\end{Not}


In application only finite data are available why we introduce additionally the following observation model, which is more realistic and has been studied for example in \cite{MR1425958, MR2240642}:

\begin{Not}[Discretized data]\label{not:discreteModel}
Let us consider the one-sided discretization of $Y_{\delta}$:
\begin{equation}\label{eq:diskrRM}
	Q Y_{\delta} = Q Ax + \delta Q \Xi = \sum_{j=1}^n Y_{\delta, j} w_{j} \text{ with } Y_{\delta, j} (\omega)
	= \left\langle y, w_{j} \right\rangle_{\mathbb{H}_{2}} + \delta \left\langle \Xi, w_{j} \right\rangle_{\mathbb{H}_{2}} (\omega), \, \omega \in \Omega,
\end{equation}
where $Q$ denotes the projection onto the linear span of an orthonormal system $ \left\{ w_{1},..., w_n \right\}$.
\end{Not}

\begin{Rem}\hfill
\begin{compactitem}
\item We assume to have observations without repetitions.
\item (\ref{eq:diskrRM}) conforms to the well-known regression model with orthonormal design.
\item It is evident that this model leads to a supplemental error term, the discretization error, which changes the convergence rates but not the underlying convergence behaviour if we require that $n=n(\delta)$ with $\underset{\delta \rightarrow 0}{\lim} \, n(\delta) = \infty $.
\end{compactitem}
\end{Rem}

To compare and qualify different methods we need an error criterion. Most authors use the mean squared error (MSE) and so will we. It is defined as follows:

\begin{Not}[Error criterion]
Let $\Xi$ satisfy assumption \ref{ass:Xi}. We set
\[
 MSE(R, \alpha, x, \delta) := \left( \mathbb{E} \left[ \left\| T^{+}y - R_{\alpha}Y_{\delta} \right\|^{2}_{\mathbb{H}_{1}} \right] \right)^{1/2} \text{, where } y = Tx.
\]
\end{Not}

\begin{Prop}[Finiteness of the mean squared error]\label{prop:endlFehler}
Let $R_{\alpha}$, $\alpha>0$, denote a regularization operator with associated regularization filter $F_{\alpha}$ satisfying (\ref{eq:Vainikko3}). If the operator $T$ is Hilbert-Schmidt, the MSE of $R_{\alpha}$ is finite for all $x \in \mathbb{H}_1$ and $\delta > 0$.
\end{Prop}

\begin{proof}
By Parseval's identity and Fubini's theorem we get for all $y \in \mathcal{D}(T^{+})$ the so called bias-variance decomposition of the mean squared error:
\begin{equation}\label{eq:BiasVar}
	\mathbb{E} \left[ \left\| T^{+}y - R_{\alpha}Y_{\delta} \right\|^{2}_{\mathbb{H}_{1}} \right] 
	= \left\| T^{+}y - R_{\alpha}y \right\|^{2}_{\mathbb{H}_{1}} 
	+ \delta^{2} \mathbb{E} \left[ \left\| R_{\alpha} \Xi \right\|^{2}_{\mathbb{H}_{1}} \right]
\end{equation}
The first term is the squared bias, which is related to the approximation error and specifies the difference between the exact solution and the expectation value of its estimate. It is finite for all $y \in \mathcal{D}(T^{+})$ and vanishes if $\alpha \rightarrow 0$ as we have shown in proposition \ref{thm:R-pktw-T+}. The variance measures the variability of the estimate caused by the noise. Applying the singular system $\left\{ \left( s_{j}; v_{j}, u_{j} \right) \right\}_{j \in \mathbb{N}} $ of $T$ with $s_j \leq \left\| T \right\|$ we get
\begin{equation}\label{eq:endlFehler}
	\mathbb{E} \left[ \left\| R_{\alpha} \Xi \right\|^{2}_{\mathbb{H}_{1}} \right]
	= \sum_{s_{j}>0} \left| F_{\alpha}(s_{j}^{2})s_{j} \right|^{2} \mathbb{E} \left[ \Xi_{j}^{2} \right] 
	\leq \sum_{s_{j} > 0} \left| F_{\alpha}(s_{j}^{2})s_{j} \right|^{2}
	\leq \tfrac{\gamma^2}{\alpha^2} \, \left\| T \right\|_{\text{HS}}^{2}.
\end{equation}
since from $ \left\| \mathrm{Cov}_{\Xi} \right\| \leq 1 $ it follows that $\mathbb{E} \left[ \Xi_{j}^{2} \right] \leq 1$ for all coordinates $\Xi_{j} := \langle \Xi, u_j \rangle_{\mathbb{H}_2}$, $j \geq 1$.
\end{proof}

\begin{Ass}\label{ass:THS}
In the following we assume the operator $T$ to be Hilbert-Schmidt and any considered regularization filter to satisfy (\ref{eq:Vainikko3}).
\end{Ass}

\begin{Rem}
We stress that the bound in (\ref{eq:endlFehler}) does not yield optimal order. 
\end{Rem}


\subsection{Regularization of statistical inverse problems}\label{sec:RegStIP}\label{sec:kindsConv}

To define convergent statistical regularization methods we need a reasonable handling of the stochastical noise when studying the asymptotic of a regularization method for $\delta \rightarrow 0$. As crucial point we recognize that not only the realization of the observations could vary for changing noise levels but even the underlying probability distribution could alter. 

\begin{Rem}[\textbf{Main idea: Linking the noise to the asymptotic of the noise level}]\label{rem:handling}
For a chosen class of probability distributions $\mathcal{W}$ we consider the asymptotic behaviour of a regularization method $(R, \alpha)$ when the index $k \geq 1$ tends to infinity, i.e. we study
\[
 \lim_{k \rightarrow \infty} \left\| T^{+}y - R_{\alpha(\delta^{(k)}, Y^{(k)})} Y^{(k)} \right\|_{\mathbb{H}_{1}}  \quad \text{ for }  y \in \mathcal{D}(T^{+}),
\]
where $Y^{(k)} := Y^{(k)}(y) := y + \delta^{(k)} \Xi^{(k)} $ with $ \Xi^{(k)} \sim \mathbb{P}^{\Xi^{(k)}} \in \mathcal{W} $, $ \delta^{(k)} > 0 $ and $\underset{k \rightarrow \infty}{\lim} \delta^{(k)} = 0 $.
\end{Rem}


\begin{Ex}\hfill
\begin{compactitem}
\item Let $\mathbb{P}^{\Xi}$ be any probability distribution and $\mathcal{W} := \lbrace \mathbb{P}^{\Xi} \rbrace$, i.e. we set $\Xi^{(k)} := \Xi$ for all $k \geq 1$. The assumed distribution can be interpreted as a priori knowledge of the noise behaviour. The most popular example of this approach are observations with Gaussian white noise. 
\item By setting $\mathcal{W} := \lbrace \mathbb{P}^{\Xi}: \Xi \sim \mathbb{P}^{\Xi} \text{ centered Hilbert-space process with } \Vert \mathrm{Cov}_{\Xi}  < \infty \rbrace$
we approve arbitrary observations $Y^{(k)} := y + \delta^{(k)} \Xi^{(k)} $ where $\Xi^{(k)} $ can be any Hilbert-space process satisfying assumption \ref{ass:Xi}. Here the change to the stochastic context causes a loss of information.
\item As a compromise we could consider any subclass of $\mathcal{W}_{0}$ such as the Dirac measures or the centered normal distributions with bounded covariance.
\end{compactitem}
\end{Ex}

\begin{Rem}[Kinds of convergence]
In order to formulate the aspired definitions we still lack in a convenient kind of convergence. In consideration of definition \ref{def:Rauschpegel} there are basically three possibilities available: convergence in mean square, convergence in probability and convergence in distribution. The latter is too weak to yield usefull results but convergence in probability should suffice for a lot of cases. Nevertheless the convergence in mean square is often prefered because of its technical advantages. One should decide as the case arises.
\end{Rem}

\begin{Def}[Convergent statistical regularization]\label{def:stKonvRV}
Let $ R $ be a linear regularization scheme, $\alpha: (0,\infty) \times \mathbb{R}^{\mathbb{N}} \rightarrow (0,\infty)$ a measurable function and $\mathcal{W}$ a class of probability distributions. We set
\begin{equation}\label{eq:Mw(y)}
	\mathbb{M}_{\mathcal{W}} (y) := \left\{ Y_{\delta} = y + \delta \Xi : \delta > 0, \, \mathbb{P}^{\Xi} \in \mathcal{W} \text{ and } \left\| \mathrm{Cov}_{\Xi} \right\| \leq 1 \right\} \text{ for any } y \in \mathcal{D}(T^{+}).
\end{equation}
The pair $\left( R, \alpha \right) $ is called convergent statistical regularization w.r.t. $\mathcal{W}$ if for all $y \in \mathcal{D}(T^{+})$ and arbitrary observations $\left\{ Y_{\delta} \right\}_{\delta > 0} \subseteq \mathbb{M}_{\mathcal{W}} (y) $ the regularized solution converges $\mathbb{P}$-stochastically to the exact one when $\delta \rightarrow 0$:
		\[
		\begin{split}
			\text{For all } \left\{ Y^{(k)} \right\}_{k \geq 1} \subseteq \mathbb{M}_{\mathcal{W}} (y) \text{ with } Y^{(k)} := y + \delta^{(k)} \Xi^{(k)} \text{ and } \lim_{k \rightarrow \infty} \delta^{(k)} = 0 \text{ we have} \\
			\lim_{k \rightarrow \infty} \mathbb{P} \left( \left\{ \omega \in \Omega: \left\| T^{+}y - R_{\alpha(\delta^{(k)}, Y^{(k)}(\omega))} Y^{(k)}(\omega) \right\|_{\mathbb{H}_{1}} > \epsilon \right\} \right) = 0  \text{ for all } \epsilon > 0.
		\end{split}
		\]
\end{Def}

\begin{Rem}
The convergence in probability could be replaced by the convergence in mean square. We call such schemes convergent statistical regularization in mean square w.r.t. $\mathcal{W}$.
\end{Rem}

\begin{Ex}\label{rem:randVar}\label{ex:konvRVenfWR}\hfill
\begin{compactitem}
\item Random variables:
Hofinger and Pikkarainen study in \cite{hofinger, MR2503292} convergence rates of the Tikhonov regularization using the Ky-Fan metric as error criterion and allowing only observations whose noise can be modeled as random variable. 
\item Statistical parameter choices:
In addition to modifications of classic parameter choices, several strategies have been developed especially for the stochastic context. One of them was introduced by Lepski{\u\i} in \cite{MR1091202} and since then adapted to various models as for example statistical inverse problems with Gaussian white noise \cite{MR2175028, MR2240642}. 
Another common parameter choice is cross-validation. In Tsybakov \cite{MR2013911} it is presented in a regression model and in \cite{MR1045442} one can find a $\delta$-free version.
\item Gaussian white noise in the abstract model (\ref{eq:ObsModel}):
In \cite{MR2175028} the convergence in mean square of a Lepski{\u\i}-type parameter choice applied to  spectral cut-off is proven for observations with white noise.
\item Gaussian white noise in the regression modell (\ref{eq:diskrRM}):
Math{\'e} and Pereverzev have shown in \cite{MR2240642} that Lepski{\u\i}'s procedure converges also with Tikhonov regularization. Our analysis in section \ref{sec:stBak} will be based on this study. That is why we want to outline briefly the crucial results. In \cite{MR2240642} the authors focused on discretized data with random noise as described in notation \ref{not:discreteModel}. They assumed that:
\renewcommand{\labelenumi}{\alph{enumi})}
\begin{compactenum}
	\item $\left\langle \Xi, w_{j} \right\rangle_{\mathbb{H}_{2}} \overset{\mathrm{iid}}{\sim} \mathcal{N}(0,1), j = 1,...,n$ 
	\item $x^+ \in T_{\varphi} (R) := \lbrace x \in \mathbb{H}_{1}: x = \varphi (T^* T) v, \Vert v \Vert \leq R \rbrace$, where $\varphi: \left( 0, \Vert T \Vert^2 \right] \rightarrow \mathbb{R}_+$ is an increasing and operator monotone function with $ \varphi (0) = 0 $.
	\item The singular values of $T$ satisfy $s_j \asymp j^{-r}$ for all $ j \geq 1$ and some $ r > 0 $.
	\item There is a constant $C > 1$ such that $ \Vert (I - Q) T: \mathbb{H}_1 \rightarrow \mathbb{H}_2 \Vert \leq C \, \mathrm{rank}(Q)^{-r} $.
\end{compactenum}
Further, they set
\renewcommand{\labelenumi}{\arabic{enumi})}
\begin{compactenum}
	\item $R_{\alpha} := ( \alpha I + B^* B )^{-1} B^*$ with $B := QT$
	\item $\alpha_0 := \delta^2$ and $ \alpha_j := \alpha_0 \, q^j, $ where $ q > 1 $ and $ j = 1,...,m := \lceil 2 \log_q (\Vert T \Vert^2 / \delta) \rceil $
	\item $x_{j, \delta} := R_{\alpha_j} Q \, y_{\delta} = ( \alpha_j I + B^* B )^{-1} B^* Q Y_{\delta} (\omega)$
	\item $n = n(\alpha) \asymp \lceil \alpha^{-1/2r} \rceil $ and $ Q = Q_{n}, $ where $\alpha > 0$ and $Q $ the described orthonormal projection onto $\mathrm{span} ( \left\lbrace w_1,...,w_n \right\rbrace )$ 
	\item Let $C_{\Psi}, \, C_1$ and $ C_2 > 0 $ be such that
	\begin{align*}
	\Psi (j) &:= C_{\Psi} \sqrt{\tfrac{1}{4 \alpha_j} \mathrm{rank} (Q)} \geq \mathbb{E} \left[ \Vert R_{\alpha_j} Q \Xi \Vert^2 \right] \text{ (decreasing) and }\\
	\Phi(j) &:= C_1 \varphi (C_2 \alpha_j) \geq \Vert T^+ y - R_{\alpha_j} Q Tx \Vert \text{ (increasing) with } j=0,...,m
	\end{align*}
	satisfy $\delta \Psi (\alpha_0) \geq \Phi (\alpha_0)$.
\end{compactenum}
Now, the regularization parameter is chosen according to 
\begin{equation}\label{eq:LepskiPW}
	\alpha_* := \alpha_{j_*} \text{ with } j_* := \max \left\lbrace j = 1,...,m: \Vert x_{k, \delta} - x_{j, \delta} \Vert \leq 4 \,  \kappa \, \delta \; \Psi (k) \text{ for all } k \leq j \right\rbrace, 
\end{equation}
where $\kappa := \sqrt{m}$. The idea of this choice is to approximate the parameter $\alpha_{\mathrm{opt}}$ which satisfies $ \delta \Psi (\alpha_{\mathrm{opt}}) = \Phi (\alpha_{\mathrm{opt}}) $. Finally, we get with $ \Theta (t) := t^{(2r+1)/4r} \varphi(t), 0 < t \leq \Vert T \Vert^2, $ and $\delta_0 > 0$ sufficiently small, that
\[
	\sup_{x \in T_{\varphi}(R)} \left( \mathbb{E} \left[ \Vert x - x_{j_*, \delta} \Vert^2 \right] \right) ^{1/2} \leq C_{\mathrm{all}} \sqrt{\lceil 2 \log_q \left( \tfrac{\Vert T \Vert^2}{\delta} \right) \rceil } \varphi \left( \Theta^{-1} \left( \tfrac{\delta}{R}\right) \right)  , \quad\delta \leq \delta_0,
\]
what converges to zero when $\delta \rightarrow 0$.
\end{compactitem}
\end{Ex}

\begin{Rem}
For more details about the concepts of \textit{general source conditions} and \textit{operator monotone functions} we refer to \cite{MR1477662, MR2384768, MR2240642} and the references therein.
\end{Rem}


\subsection{Relation between classic and statistical regularization methods}\label{sec:transf}

As justification for section \ref{sec:RegStIP} and as preparation of section \ref{sec:stBak} we are interested in the connection of regularization methods of the two settings. 
In general, we have to modify at least the parameter choice $\alpha$ because of the changed domain of definition. In order to formulate sufficient criteria for the stochastical convergence of $(R, \tilde{\alpha})$ we need to control the decay of $\tilde{\alpha}( \delta, Y_{\delta})$. The following notation will help us to describe it conveniently.

\begin{Not}[Stochastic Landau-Symbol ]\label{not:Landau}
Let $(Z_{n})_{n \in \mathbb{N}}$ be a sequence of random variables on a probability space $(\Omega, \mathcal{F},\mathbb{P})$ and $(c_{n})_{n \in \mathbb{N}}$ a sequence of real-valued constants. We denote
\[
	Z_{n} = o_{\mathbb{P}}(c_{n}) :\Longleftrightarrow \lim_{n \rightarrow \infty} \mathbb{P} \left( \vert \tfrac{Z_{n}}{c_{n}} \vert > \epsilon \right)  = 0 \text{ for all } \epsilon > 0.
\]
\end{Not}


\begin{Prop}\label{thm:det-stRV}
Let  $(R, \alpha)$ be any generally convergent regularization,
\begin{equation}\label{eq:W0}
	\mathcal{W} \subseteq \lbrace \mathbb{P}^{\Xi}: \Xi \sim \mathbb{P}^{\Xi} \text{ centered Hilbert-space process with} \Vert \mathrm{Cov}_{\Xi}  \Vert < \infty \rbrace =: \mathcal{W}_0
\end{equation}
 and $\mathbb{M}_{\mathcal{W}} (y) $, $ y \in \mathcal{D}(T^{+}) $, such as in (\ref{eq:Mw(y)}). The modified method $(R, \tilde{\alpha})$ constitutes a convergent statistical regularization w.r.t. $\mathcal{W}$ for any measurable function $\tilde{\alpha}: (0,\infty) \times \mathbb{R}^{\mathbb{N}} \rightarrow (0, \infty)$ if for arbitrary observations $\left\{ Y^{(k)} \right\}_{k \geq 1} $ with $Y^{(k)} := y + \delta^{(k)} \Xi^{(k)} \in \mathbb{M}_{\mathcal{W}} (y)$ it holds
\begin{equation}\label{eq:a=oPd}
  \lim_{k \rightarrow \infty} \mathbb{P} \left( \tilde{\alpha} ( \delta_k, Y_k) > \epsilon \right) = 0 \text{ for all } \epsilon > 0 \quad \text{ and } \quad (\tilde{\alpha} ( \delta_k, Y_k))^{-1} = o_{\mathbb{P}} (\delta_k^{-1}).
\end{equation}
\end{Prop}

\begin{proof}
Let $ y \in \mathcal{D}(T^{+}) $, $\left\{ Y^{(k)} \right\}_{k \geq 1} \subseteq \mathbb{M}_{\mathcal{W}} (y)$ with $Y^{(k)} := y + \delta^{(k)} \Xi^{(k)} $ and $\underset{k \rightarrow \infty}{\lim} \delta^{(k)} = 0$. Proposition \ref{prop:endlFehler} yields with assumption \ref{ass:THS} for any number $\alpha>0$ the finiteness of the mean squared error:
\[
	\mathbb{E} \left[ \left\| T^{+}y - R_{\alpha} Y^{(k)} \right\|^{2}_{\mathbb{H}_{1}} \right]
	\leq \left\| T^{+}y - R_{\alpha} y \right\|^{2}_{\mathbb{H}_{1}} + (\delta^{(k)})^{2} \tfrac{\gamma^2}{\alpha^2} \, \left\| T \right\|_{\text{HS}}^{2} < \infty\, , \quad k \geq 1.
\]
Now, we consider a measurable function $\tilde{\alpha}: (0,\infty) \times \mathbb{R}^{\mathbb{N}} \rightarrow (0,\infty)$ satisfying (\ref{eq:a=oPd}) and insert in place of the number $\alpha$ the function value $ \tilde{\alpha}(\delta^{(k)}, Y^{(k)}(\omega)) $, where $Y^{(k)} (\omega) = \lbrace Y^{(k)}_{\delta^{(k)}, j} (\omega) \rbrace_{j \geq 1}$ for $\omega \in \Omega$ and $k \geq 1$. In doing so we allow for a moment that the parameter choice and the regularization operator are applied to different realizations of $Y^{(k)}$, $k \geq 1$. We get from proposition \ref{thm:R-pktw-T+} that
\begin{equation}\label{eq:stKonvMSE}
	\lim_{k \rightarrow \infty} \mathbb{P} \left( \left\{ \omega \in \Omega: \mathbb{E} \left[ \left\| T^{+}y - R_{\tilde{\alpha}(\delta^{(k)}, Y^{(k)}(\omega))} Y^{(k)} \right\|^{2}_{\mathbb{H}_{1}} \right] > \epsilon \right\} \right) = 0 
\end{equation}
 for any $\epsilon > 0$ since the sum of two stochastical convergent sequences converges stochastically. 
So, we can say: For all $\epsilon > 0$ there exists a subset $\Omega_{\epsilon} \subseteq \Omega$ with $\mathbb{P}(\Omega_{\epsilon}) \geq 1 - \epsilon$, such that 
\[
	\lim_{k \rightarrow \infty} \mathbb{E} \left[ \left\| T^{+}y - R_{\tilde{\alpha}(\delta^{(k)}, Y^{(k)}(\tilde{\omega}))} Y^{(k)} \right\|^{2}_{\mathbb{H}_{1}} \right] = 0
	\quad \text{for all} \, \tilde{\omega} \in \Omega_{\epsilon} \text{ with } \mathbb{P}({\tilde{\omega}}) > 0.
\]
Further we can deduce that for all $\epsilon,\eta > 0$ and $ \tilde{\omega} \in \Omega_{\epsilon} $ with $\mathbb{P}({\tilde{\omega}}) > 0$
\[
	\lim_{k \rightarrow \infty} \mathbb{P} \left( \left\{ \omega \in \Omega: \left\| T^{+}y - R_{\tilde{\alpha}(\delta^{(k)}, Y^{(k)}(\tilde{\omega}))} Y^{(k)}(\omega) \right\|^{2}_{\mathbb{H}_{1}} > \eta \right\} \right) = 0.
\]
Finally, we achieve
\[
	\lim_{k \rightarrow \infty} \mathbb{P} \left( \left\{ \omega \in \Omega: \left\| T^{+}y - R_{\tilde{\alpha}(\delta^{(k)}, Y^{(k)}(\omega))} Y^{(k)}(\omega) \right\|^{2}_{\mathbb{H}_{1}}> \eta \right\} \right) \leq \epsilon
\]
for all $\epsilon, \eta > 0$. Since $\epsilon$ is independent of $\eta$, we can conclude stochastical convergence.
\end{proof}

\begin{Prop}\label{prop:TransfRVrv} 
Any generally convergent regularization $(R, \alpha)$, where $\alpha$ is measurable, satisfies definition \ref{def:stKonvRV} of convergent statistical regularization w.r.t. 
\begin{equation}\label{eq:W2}
	\mathcal{W}_{1} \subseteq \mathcal{W}_{2} := \lbrace \mathbb{P}^{\Xi}: \Xi \in L^{2}(\Omega, \mathbb{H}_{2}) \text{ with } \mathbb{E} \left[ \left\| \Xi \right\|_{\mathbb{H}_{2}}^{2} \right] \leq 1 \rbrace.
\end{equation}
The converse holds if $\mathcal{W}_{1}$ contains the Dirac measures.
\end{Prop}


\begin{proof}
Let $(R, \alpha)$ be a generally convergent regularization method with measurable $\alpha$, $ y \in \mathcal{D}(T^{+}) $ and $\left\{ Y^{(k)} \right\}_{k \geq 1} $ with $Y^{(k)} := y + \delta^{(k)} \Xi^{(k)} $, $\mathbb{P}^{\Xi^{(k)}} \in \mathcal{W}_{1}$, $\delta^{(k)} > 0$ for all $k \geq 1$ and $\underset{k \rightarrow \infty}{\lim} \delta^{(k)} = 0$. 
We fix $\epsilon > 0$, set $C := \epsilon^{-1/2}$ and define for any $\omega \in \Omega$ the set
\[
\begin{split}
	\mathcal{K}(\omega) &:= \left\{ k \geq 1: \left\| \Xi^{(k)}(\omega) \right\|_{\mathbb{H}_{2}} \leq C \right\} \subseteq \mathbb{N} \text{ and the number }\\
	k(\omega) &:= \text{argmin} \left\{ k \geq 1: \left\| T^{+}y - R_{\alpha(\delta^{(k)}, Y^{(k)}(\omega))} Y^{(k)}(\omega) \right\|_{\mathbb{H}_{1}} \leq \epsilon \; \forall \, l \in \mathcal{K}(\omega) \text{ with } l \geq k \right\}.
\end{split}
\]
Then it follows from Chebychev's inequality and the convergence of $(R, \alpha)$ that
\[
\begin{split}
	\mathbb{P} &\left( \left\{ \omega \in \Omega: \left\| T^{+}y - R_{\alpha(\delta^{(k)}, Y^{(k)}(\omega))} Y^{(k)}(\omega) \right\|_{\mathbb{H}_{1}} > \epsilon \right\} \right)\\
	&\quad \leq \mathbb{P} \left( \left\{ \omega \in \Omega: k \notin \mathcal{K}(\omega) \right\} \right)
	+ \mathbb{P} \left( \left\{ \omega \in \Omega: k < k(\omega) \right\} \right) < 2 \epsilon
\end{split}
\]
 for $k \geq 1$ sufficiently large and finally
\[
	\lim_{k \rightarrow \infty} \mathbb{P} \left( \left\{ \omega \in \Omega: \left\| T^{+}y - R_{\alpha(\delta^{(k)}, Y^{(k)}(\omega))} Y^{(k)}(\omega) \right\|_{\mathbb{H}_{1}} > \epsilon \right\} \right) = 0  \text{ for all } \epsilon > 0.
\]
\end{proof}

\begin{Prop}\label{prop:st-detRVnf}
Any purely data driven convergent statistical regularization $(R, \alpha)$ w.r.t. $\mathcal{W}_{0}$, induces a purely data driven generally convergent regularization $(R, \tilde{\alpha})$.
\end{Prop}

\begin{proof}[Proof sketch]
Let us contemplate deterministic observations of the form $y^{(k)} := y + \delta^{(k)} \xi^{(k)} \in \mathbb{H}_{2}$ with $y \in \mathcal{D}(T^{+})$, $\left\| \xi^{(k)} \right\|_{\mathbb{H}_{2}} \leq 1$, $\delta^{(k)} > 0$ for all $k \geq 1$ and $\underset{k \rightarrow \infty}{\lim} \delta^{(k)} = 0$. We define for any $k \geq 1$ the following Hilbert-space valued random variable
\[
	Y^{(k)}(\omega) := 
	\begin{cases}
		y^{(k)}, &\text{ if } \omega \in \Omega_{1} \\
		-y^{(k)}, &\text{ if } \omega \in \Omega_{2},
	\end{cases}
\]
where $\mathbb{P}(\Omega_{1}) = \mathbb{P}(\Omega_{2}) = 0.5$. Every random variable $ Y^{(k)}$, $k \geq 1 $, can be identified with a centered Hilbert-space process, such that the function
\[
	\tilde{\alpha}: \mathbb{H}_{2} \rightarrow (0, \infty), \, y^{(k)} \mapsto \alpha \left( \left\{ y_{\delta^{(k)}, j}^{(k)} \right\}_{j \geq 1} \right),
\]
where $y_{\delta^{(k)}, j}^{(k)} := Y_{\delta^{(k)}, j}^{(k)} (\omega)$ for any $\omega \in \Omega_{1}$ and $ j \geq 1$, constitutes with the regularization $R$ a purely data driven generally convergent regularization.
\end{proof}

\begin{Rem}\label{rem:GenerSt-detRVnf}
The proposition holds also for methods w.r.t. a subclass $ \mathcal{W} \subseteq \mathcal{W}_{0}$ if $\mathcal{W}$ allows for arbitrary deterministic observations $\left\lbrace y^{(k)} \right\rbrace_{k \geq 1} $ of the above form the definition of a sequence $\left\lbrace Y^{(k)} \right\rbrace_{k \geq 1} \subseteq \mathbb{M}_{\mathcal{W}} (y)$ with $ \mathbb{P} \left( \left\lbrace \omega \in \Omega: Y_k (\omega) = y_k \right\rbrace \right) > \eta$ for $ \eta>0$.
\end{Rem}


\section[The Bakushinskii veto for statistical inverse problems]{The Bakushinski{\u\i} veto for statistical inverse problems}	\label{sec:stBak}\label{sec:RegUnabhRP}

The following study was motivated by the paper \"{}Regularization independent of the noise level: an analysis of quasi-optimality\"{} by Bauer and Rei{\ss} \cite{MR2438944}, which raised the question of the transferability of the Bakushinski{\u\i} veto to statistical inverse problems.

\begin{Thm}\label{thm:stBak}\label{thm:stBakfPXiFix}\hfill
\begin{compactenum}
 \item A purely data driven convergent statistical regularization method w.r.t. $\mathcal{W}_{0} $, see (\ref{eq:W0}), exists if and only if the range $\mathcal{R}(T)$ of the operator $T$ is closed.
 \item For certain probability distributions $\mathbb{P}^{\Xi}$ there exist purely data driven convergent statistical regularization w.r.t. $\mathcal{W} := \lbrace \mathbb{P}^{\Xi} \rbrace$ even if the problem is ill-posed.
\end{compactenum}
\end{Thm}

\begin{Rem}[Generalization]
The first statement remains valid for sufficiently large subclasses of $\mathcal{W}_{0}$ such as $\mathcal{W}_{2}$ of (\ref{eq:W2}) or the class of all Dirac measures. We refer to remark \ref{rem:GenerSt-detRVnf}.
\end{Rem}

For the proof of the second statement we need some preperation:

\begin{Not}[Setting]\label{not:setting}
In order to construct an example supporting theorem \ref{thm:stBakfPXiFix} (2) let us focus on an operator $T: L^2( \left[ 0, 1 \right] ) \rightarrow L^2( \left[ 0, 1 \right] )$ and data with Gaussian white noise modeled by
$Y_{\delta} (t) = Tx(t) + \delta \, \Xi_{t}, \, t \in \left[ 0,1 \right] ,$
which is consistent with (\ref{eq:ObsModel}). 
We consider the equidistant decomposition $\mathcal{Z}_n := (0=t_0 < t_1 < ... < t_n = 1) \text{ with } t_j := j / n \text{ for } j= 0,...,n$ and the orthornormal system $ \lbrace \varphi_j \rbrace_{j=1,...,n}$, where 
$ \varphi_j := \left( \sqrt{t_j - t_{j-1}} \right)^{-1} \chi_{\left[ t_{j-1}, t_j \right)}. $
By projecting $ Y_{\delta} $ onto the linear span of $ \lbrace \varphi_j \rbrace_{j=1,...,n}$ we get a finite set of coefficients
\[
	Y_{\delta,j} := \langle Y_{\delta}, \varphi_j \rangle_{L^2( \left[ 0, 1 \right] )}
	= \left( \sqrt{t_j - t_{j-1}} \right)^{-1} \, \int_{t_{j-1}}^{t_j} Tx(s)ds + \delta \epsilon_j , \quad \epsilon_j \overset{iid}{\sim} \mathcal{N}(0,1)
\]
with $j= 1,...,n$, such that
\[
	QY_{\delta} (t) = \sum_{j = 1}^{n} Y_{\delta,j} \varphi_j (t) = (t_j - t_{j-1})^{-1} \, \int_{t_{j-1}}^{t_j} Tx(s)ds + \sqrt{n} \, \delta \, \epsilon_j \text{ for } t \in \left[ t_{j-1}, t_j \right).
\]
\end{Not}

\begin{Rem}[Outline]
This setting conforms to the regression model with orthonormal design and without repetitions as discribed in notation \ref{not:discreteModel}. In example \ref{ex:konvRVenfWR} we mentioned that Tikhonov regularization forms with a Lepski{\u\i}-type parameter choice a convergent statistical regularization method $(R, \alpha)$ w.r.t. $\mathcal{N}(0,I)$ \cite{MR2240642}. Plugging in an estimation of the noise level into this method we can deduce a purely data driven one as we will verify now.
\end{Rem}

For that purpose we want to use the following estimator:
\begin{Def}[The estimator \cite{MR897854}]\label{def:estimator}
\begin{equation}\label{eq:SchaetzeRPimRegrM}
	\tilde{\delta}_n^{2} := \frac{1}{2n^2} \sum_{j = 1}^{n}  ( Q Y_{\delta}(t_{j}) - Q Y_{\delta}(t_{j-1}) )^{2}
\end{equation}
\end{Def}

Before studying its asymptotical behaviour we remind of the following notation:

\begin{Not}(\cite{MR1902050})
Let $I$ denote an interval. A function $f:I \rightarrow \mathbb{R}$ is called H\"{o}lder continuous with exponent $s \in \left( 0,1 \right] $ if for all $t_{0} \in I$ a neighborhood $U \subseteq \mathbb{R}$ exists, such that 
\[
	\sup_{t, t' \in U \cap I, t \neq t'} \dfrac{\vert f(t) - f(t') \vert }{ \vert t - t' \vert^{s} } < \infty.
\]
\end{Not}

\begin{Ass}\label{ass:TxHcont}
Let $y = Tx \in L^2( \left[ 0, 1 \right] )$ be H\"{o}lder-continuous of order $s \in \left( \tfrac{1}{2}, 1 1 \right]$.
\end{Ass}

\begin{Ex}
Assumption \ref{ass:TxHcont} is satisfied for any integral operator $T: L^2( \left[ 0, 1 \right] ) \rightarrow L^2( \left[ 0, 1 \right] )$,
\[
	(Tx)(t) = \int_0^t k(t,u) x(u) du,
\]
 with kernel $k: \left[ 0, 1 \right]^2 \rightarrow \mathbb{R}$ satisfying for some constant $C > 0$
\[
	\sup_{u \in \left[ 0, 1 \right]} \vert k(t,u) - k(t',u) \vert \leq C \vert t - t' \vert^s, \, t,t' \in \left[ 0, 1 \right] .
\]
\end{Ex}

\begin{Conc}
Assumption \ref{ass:TxHcont} implies in pursuance of \cite[pages 212-213]{MR606200} that 
\[
	( Q y (t_{j-1}) - Q y (t_{j}) )^{2} \asymp O (n^{-2s}), \quad j=1,...,n,
\]
what from we can deduce the asymptotical unbiasedness of $\tilde{\delta}_n^{2}$ when $n \rightarrow \infty$:
	$$\mathbb{E} (\tilde{\delta}^{2}_{n}) 
	= \delta^2 + O(n^{-(1+2s)}) \text{ and } \mathbb{E} (\tilde{\delta}^{2}_{n}) \geq \delta^2 $$
\end{Conc}

\begin{Rem}
In proposition \ref{prop:MPmitNneu} we need $s > \tfrac{1}{2}$.
\end{Rem}

\begin{Prop}[Concentration inequality]\label{prop:Omega-}
Let $n = n(\delta) \asymp \lceil \delta^{-\eta} \rceil$ with $2 > \eta \geq \tfrac{2}{1+2s}$, $ \hat{\delta} := \tau \, \tilde{\delta}_n $ and
\begin{equation}\label{eq:Omega+}
	\Omega_+ := \Omega_+ (\delta, \tau, K) := \left\lbrace \omega \in \Omega: \hat{\delta} (\omega) \in \left[ \delta, K \tau \delta \right] \right\rbrace, \, \tau, K >1 \text{ appropriate.}
\end{equation}
The following assertions hold for all $\delta \leq \delta_0$ with $\delta_0>0$ sufficiently small:
\begin{compactenum}
\item There are constants $C_1, C_2 > 0$ such that
	$
		\mathbb{P} ( \Omega \setminus \Omega_+ ) 
		\leq C_1 \exp \left( - \, C_2 \, n^2(\alpha, \delta) \,\delta^2 \right).
	$
\item It holds for $\alpha > 0$ and some $C_3>0$ that
	\[
		\underset{T^+ y \in T_{\varphi}(R)}{\sup} \int_{\Omega} \Vert T^+ y - R_{\alpha} Q Y_{\delta} \Vert^4 d \mathbb{P}
		\leq C_3 \, \delta^{4} \alpha^{-4} n^2(\alpha, \delta).
	\]
\end{compactenum}
\end{Prop}

We want to use the following Lemma for the proof of proposition \ref{prop:Omega-}:

\begin{Lem}\label{lem:equivMoments}\hfill
\begin{compactenum}
\item Let $X$ be a Gaussian random vector in a Banach space $\mathcal{B}$ and $\Vert X \Vert_p := \mathbb{E} \left[ \Vert X \Vert_{\mathcal{B}}^p \right]^{1/p} $, $0<p<\infty$, the $L^p$-norm of $X$. For all $0 < p,q < \infty$ there is a constant $K_{p,q} > 0$ such that $ \Vert X \Vert_p \leq K_{p,q} \Vert X \Vert_q $.
\item Let $X \in L^4(\Omega, \mathcal{F}, \mathbb{P})$ be nonnegative. It holds
\[
	\int_{\Omega} X^4(\omega) d \mathbb{P} (\omega)
	= 4 \int_0^{\infty} t^3 \, \mathbb{P} (X > t) dt.
\]
\item $ \Vert (\alpha I + B^* B )^{-1} B^* B \Vert \leq 1$
\end{compactenum}
\end{Lem}

\begin{proof}[Proof  of lemma \ref{lem:equivMoments}]
The first statement can be found in \cite[page 60]{MR1102015} and the second one follows by a generalized version of partial integration, which is given in \cite[chapter 5, $\S$ \medskip 6]{MR0270403}. By spectral calculus we get the inequality in (3).
\end{proof}

\begin{proof}[Proof  of proposition \ref{prop:Omega-}]\hfill
\begin{compactenum}
\item Using Lemma \ref{lem:equivMoments} (1) we get with $\tilde{\delta}_n = \Vert X_n \Vert $, where
	\[
		X_n := \tfrac{1}{\sqrt{2} \, n} \left( Q Y_{\delta} (t_1) - Q Y_{\delta} (t_0), ..., Q Y_{\delta} (t_{n}) - Q Y_{\delta} (t_{n-1}) \right) \sim \mathcal{N} (\mathbb{E} X_n, I_n) ,
	\]
	that
	\[
	\begin{split}
		\mathbb{P} &( \Omega \setminus \Omega_+ ) 
		= \mathbb{P} \left( \left\lbrace \omega \in \Omega: \tau \, \tilde{\delta}_n (\omega) < \delta \right\rbrace \right) 
		+ \mathbb{P} \left( \left\lbrace \omega \in \Omega: \tau \, \tilde{\delta}_n  (\omega) > \tau \, K \delta \right\rbrace \right) \\
		&= \, \mathbb{P} \left( K_{2,1}^{-1} \Vert \tilde{\delta}_n \Vert_{2} - \tilde{\delta}_n > K_{2,1}^{-1} \Vert  \tilde{\delta}_n \Vert_{2} - \tfrac{\delta}{\tau} \right) 
		+ \,\mathbb{P} \left( \tilde{\delta}_n - K_{1,2} \Vert  \tilde{\delta}_n \Vert_{2} > K \delta - K_{1,2} \Vert \tilde{\delta}_n \Vert_{2} \right) \\
		&\leq \, \mathbb{P} \left( \vert \mathbb{E} \left[ \tilde{\delta}_n \right] - \tilde{\delta}_n \vert > K_{2,1}^{-1} \Vert  \tilde{\delta}_n \Vert_{2} - \tfrac{\delta}{\tau} \right) 
		+ \, \mathbb{P} \left( \vert \tilde{\delta}_n - \mathbb{E} \left[ \tilde{\delta}_n \right] \vert > K \delta - \Vert  \tilde{\delta}_n \Vert_{2} \right)
	\end{split}
	\]
	since the Cauchy-Schwarz inequality yields $K_{1,2} =1$. At this point, we would like to apply the concentration inequality (3.2) in \cite[page 57]{MR1102015} what for we have to ensure that $\tau \Vert \tilde{\delta}_n \Vert_{2} > K_{2,1}\delta$ and $K \delta > \Vert  \tilde{\delta}_n \Vert_{2}$.
	The first requirement is satisfied for all $\tau > K_{2,1}$ as $\Vert \tilde{\delta}_n \Vert_{2} \geq \delta$. For the second we need that $K \gg 1$ since we have for a constant $c > 0$
	\[
		\Vert \tilde{\delta}_n \Vert_{2} \leq \delta + c n^{-(1+2s)/2} \asymp \delta + c \, \delta^{\eta (1+2s)/2} \leq (c+1) \, \delta \text{ for } \delta \in (0,1).
	\]
	Supposing that $\tau$ and $K$ are appropriate it follows that for some constants $C_1, C_2 > 0 $
	\[
	\begin{split}
		\mathbb{P} &( \Omega \setminus \Omega_+ ) \\
		&\leq \, 2 \exp \left(  - \tfrac{2}{\pi^2} n^2(\alpha, \delta) \left[  K_{2,1}^{-1} \Vert  \tilde{\delta}_n \Vert_{2} - \tfrac{\delta}{\tau} \right]^2 \right) 
		+ \, 2 \exp \left(  - \tfrac{2}{\pi^2} n^2(\alpha, \delta) \left[  K \delta - \Vert \tilde{\delta}_n \Vert_{2} \right]^2 \right)\\
		&\leq C_1 \exp \left( - \, C_2 \, n^2(\alpha, \delta) \,\delta^2 \right).
	\end{split}
	\]
\item After lemma \ref{lem:equivMoments} (2) and (3) and the concentration inequality (3.5.) in \cite[page 59]{MR1102015} it holds for some constant $ C_3 > 0 $ that
	\[
	\begin{split}
		\int_{\Omega} &\Vert T^+ y - R_{\alpha} Q Y_{\delta} \Vert^4 d \mathbb{P}
		= 4 \int_0^{\infty} t^3 \mathbb{P} \left( \Vert T^+ y - R_{\alpha} Q Y_{\delta} \Vert > t \right) dt \\
		&\leq 4 \int_0^{2 \Vert T^+ y \Vert_{\mathbb{H}_1}} t^3 dt 
		+ 16 \int_{2 \Vert T^+ y \Vert_{\mathbb{H}_1}}^{\infty} t^3 \, \mathbb{P} \left( 2 \Vert T^+ y \Vert + \delta \Vert R_{\alpha} Q \Xi \Vert > t \right)  dt \\
		&\leq 16 \Vert T^+ y \Vert_{\mathbb{H}_1}^4
		+ 16 \int_{2 \Vert T^+ y \Vert_{\mathbb{H}_1}}^{\infty} t^3 \exp \left( - \tfrac{(t - 2 \Vert T^+ y \Vert_{\mathbb{H}_1})^2}{E} \right) dt \\
 		&= 16 \Vert T^+ y \Vert_{\mathbb{H}_1}^4
 		+ 16 \int_{0}^{\infty} \left( t^3 + 6 t^2 \Vert T^+ y \Vert_{\mathbb{H}_1} + 12 t \Vert T^+ y \Vert_{\mathbb{H}_1}^2 + 8 \Vert T^+ y \Vert_{\mathbb{H}_1}^3 \right) e^{- t^2 / E} dt \\
		&= 16 \Vert T^+ y \Vert_{\mathbb{H}_1}^4 + \tfrac{1}{2} E^2 + \tfrac{ 3 \sqrt{\pi} }{2} \Vert T^+ y \Vert_{\mathbb{H}_1} E^{3/2} + 6 \Vert T^+ y \Vert_{\mathbb{H}_1}^2 E + 4 \sqrt{\pi} \Vert T^+ y \Vert_{\mathbb{H}_1}^3 \sqrt{E} \\
		&\leq 16 \Vert T^+ y \Vert_{\mathbb{H}_1}^4 + C_3 \, \delta^{4} \alpha^{-4} n^2(\alpha, \delta),
	\end{split}
	\]
	where $0 < s_j \leq \Vert T \Vert $ yields
	\[
		E := 8 \mathbb{E} \left[ \Vert \delta R_{\alpha} Q \Xi \Vert^2 \right]
		= 8 \delta^2 \sum_{s_j>0} \vert \tfrac{s_j}{\alpha + s_j^2} \vert^2 \mathbb{E} \left[ \vert \langle Q \Xi, u_j \rangle_{\mathbb{H}_2} \vert^2 \right]
		\leq 8 \delta^{2} \Vert T \Vert^2 \alpha^{-2} n(\alpha, \delta).
	\]
	The assertion follows for all $\delta \leq \delta_0$ with $\delta_0 > 0$ sufficiently small.
\end{compactenum}
\end{proof}

\begin{Rem}[Asymptotic behaviour of $n$]
In example \ref{ex:konvRVenfWR} we set $n \asymp \lceil \alpha^{-1/2r} \rceil$, $r>0$, but in proposition \ref{prop:Omega-} it was more advantageous to link $n$ to $\delta$. Combining the two approaches we get
\begin{equation}\label{eq:n(a,d)}
	n := n(\alpha, \delta) = \max \lbrace n_1(\alpha), n_2(\delta) \rbrace \text{, where } 
	n_1(\alpha)  \asymp \lceil \alpha^{-1/2r} \rceil \text{ and } 
	n_2(\delta) \asymp \lceil \delta^{-\eta} \rceil.
\end{equation}
\end{Rem}

Due to the fact that we take another asymptotic behaviour of $n$ as basis of our analysis than stated in example \ref{ex:konvRVenfWR} we have to revise the convergence result.

\begin{Prop}\label{prop:MPmitNneu}
Let $\alpha_* := \alpha_{j_*} (\delta, Y_{\delta})$ denote the regularization parameter according to Lepski{\u\i}'s principle as described in example \ref{ex:konvRVenfWR}. If we assume that $n = n(\alpha, \delta)$ as in (\ref{eq:n(a,d)}) with $\eta < 2$ (instead of $n = n(\alpha) \asymp \lceil \alpha^{-1/2r} \rceil$ as before) then
\[
	\lim_{\delta \rightarrow 0} \left( \sup_{T^+ y \in T_{\varphi}(R)} \mathbb{E} \left[ \Vert T^+ y - R_{\alpha_*} Q Y_{\delta} \Vert^2 \right] \right) = 0.
\]
\end{Prop}

\begin{proof}
Math\'{e} and Pereverzev have shown in \cite[Theorem 5]{MR2240642}  that under the assumptions and notations of example \ref{ex:konvRVenfWR} it holds for some $C_0 > 0$ that
\[
	\sup_{T^+ y \in T_{\varphi}(R)} \mathbb{E} \left[ \Vert T^+ y - R_{\alpha_*} Q Y_{\delta} \Vert^2 \right] 
	\leq C_{0} \sqrt{ \lceil 2 \log_q (\Vert T \Vert^2/ \delta) \rceil } \varphi (\check{\alpha}),
\]
where $\check{\alpha} := \alpha_{\check{j}}$ with $\check{j} := \max \left\lbrace j \leq m: \Phi(j) \leq \delta \Psi(j) \right\rbrace $. The proof of this bound does not depent on the asymptotic behaviour of $n$ aside from the requirement of the existence of a constant $D>0$ satisfying $\Psi(j) \leq D \, \Psi (j+1)$ for all $j=0,...,m-1$. Since this is fulfilled even for our new choice of $n$ we cite the given inequality without further proof. The only modification which we made is a slight change of the definition of $\check{j}$, which simplifies the notation. Now, we want to prove that the right hand side converges to zero. We follow the ideas in \cite{MR2240642} and set  
\[
	\Theta_{\delta} (t) := \max \lbrace \lceil t^{1/4s} \rceil, \lceil \delta^{\eta/2} \rceil \rbrace \sqrt{t} \varphi(t), \, t>0, \text{ and } 
	\alpha^*:= \inf \left\lbrace \alpha > 0: \Theta_{\delta} (\alpha) \geq \delta \right\rbrace.
\]
$\Theta_{\delta}$ is increasing in $t$ such that for every $\delta > 0$ there is a unique choice for $\alpha^*$. We notice that
\[
	\delta \, \Psi(\alpha^*) = \delta \, C_{\Psi} \sqrt{\tfrac{n(\alpha^*, \delta)}{4 \alpha^*} } 
	= \tfrac{C_{\Psi}}{2} \, \delta \,  \left( \max \lbrace \lceil (\alpha^*)^{1/4s} \rceil, \lceil \delta^{\eta/2} \rceil \rbrace \sqrt{\alpha^*} \right)^{-1}
	\leq \tfrac{C_{\Psi}}{2} \, \varphi(\alpha^*).
\]
This leads to $\check{\alpha} \leq \alpha^{*}$ because of the definition of $\Phi$ and the monotonicity of $\Phi$ and $\Psi$. Finally, we can deduce
\begin{equation}\label{eq:MPnewRate}
	\lim_{\delta \rightarrow 0} \left( \sup_{T^+ y \in T_{\varphi}(R)} \mathbb{E} \left[ \Vert T^+ y - R_{\alpha_*} Q Y_{\delta} \Vert^2 \right] \right)
	\leq \lim_{\delta \rightarrow 0} C_{\mathrm{all}} \sqrt{ \lceil 2 \log_q (\Vert T \Vert^2/ \delta) \rceil } \varphi (\alpha^*(\delta)) = 0
\end{equation}
since $\underset{\delta \rightarrow 0}{\lim} \alpha^*(\delta) = 0$ if $\eta < 2$.
\end{proof}

\begin{Rem}
The convergence rate given in (\ref{eq:MPnewRate}) has not to be optimal.
\end{Rem}

Finally, we achieve:

\begin{proof}[Proof of theorem \ref{thm:stBakfPXiFix}]
Any purely data driven convergent statistical regularization method $(R, \alpha )$ w.r.t. $\mathcal{W}_{0}$ induces the existence of a purely data driven convergent regularization $(R, \tilde{\alpha})$ in terms of definition \ref{def:konvRV}, as shown in proposition \ref{prop:st-detRVnf}. If so, the range $\mathcal{R}(T)$ of $T$ is closed, see lemma \ref{rem:detBak}. So, we turn to the second statement:\\
We consider the setting desribed in notation \ref{not:setting} with $Tx$ satisfying assumption \ref{ass:TxHcont}, the estimator $\tilde{\delta}_n^2$ given in definition \ref{def:estimator} and the set $\Omega_+$ introduced in (\ref{eq:Omega+}). Let $R_{\alpha}$,  $m$, $\lbrace \alpha_j \rbrace_{j=0,...,m}$, $\lbrace x_{j, \delta}  \rbrace_{j=0,...,m}$, $T_{\varphi} (R)$, $\Psi$ and $\kappa$ be as in example \ref{ex:konvRVenfWR} and $n := n(\alpha, \delta) $ as in (\ref{eq:n(a,d)}). First of all we want to verify if the assumptions of example \ref{ex:konvRVenfWR} are satisfied. The first one follows by definition and the second if $x^+ \in T_{\varphi} (R)$. The definition of the projection $Q$ and the H\"{o}lder continuity of $Tx$ yield by \cite[pages 212-213]{MR606200} asumption (d) since
\[
	\Vert (I - Q) T: \mathbb{H}_1 \rightarrow \mathbb{H}_2 \Vert \leq C \, \mathrm{rank}(Q)^{-s},
\]
where $s \in (1/2,1]$ denominates the H\"{o}lder exponent of $Tx$. Assumption (c) has been used in \cite{MR2240642} as basis of assumption (d) and in order to prove the order optimality of the convergence result, why we can ignore it. As a consequence we set $r:=s$ in $n_1(\alpha) \asymp \lceil \alpha^{-1/2r} \rceil$ such that
\[
	n_1 (\alpha) \leq n_1 (\alpha_0) \asymp \delta^{-1/s} \leq \delta^{-\eta} \asymp n_2(\delta) \text{ if } \eta \geq 1/s.
\]
Now, we want to examine
\[
	\mathbb{E} \left[ \Vert T^+ y - R_{\hat{\alpha}_*} Q Y_{\delta} \Vert^2 \right] = \int_{\Omega_{+}} \Vert T^+ y - R_{\hat{\alpha}_*} Q Y_{\delta} \Vert^2 d \mathbb{P} + \int_{\Omega \setminus \Omega_+} \Vert T^+ y - R_{\hat{\alpha}_*} Q Y_{\delta} \Vert^2 d \mathbb{P},
\]
where
$
	\hat{\alpha}_* := \alpha_{j_*} (\hat{\delta}, Y_{\delta}) \text{ with } \hat{\delta} := \tau \tilde{\delta}_n
$
denotes the regularization parameter resulting from Lepski{\u\i}'s principle (\ref{eq:LepskiPW}) when using the estimated noise level. It is quite evident that 
\[
	\Omega_+ \subseteq \Omega_{\kappa} := \Omega_{\kappa} (\delta) := \left\lbrace \omega \in \Omega: \max_{j=1,...,m(\delta)} \tfrac{\delta \, \Vert R_{\alpha_j} Q_n \, \Xi (\omega) \Vert_{\mathbb{H}_2}}{\Psi_{\delta}(j)} \leq \kappa \right\rbrace
\]
if the constant $C_{\Psi} >0$ in $\Psi$ is sufficiently large. As $\alpha_{j_*} (\delta, Y_{\delta})$ and $\hat{\alpha}_*$ lead on $\Omega_+$ to the same asymptotic behaviour of $R_{\alpha} Q Y_{\delta}$ we can deduce from proposition \ref{prop:MPmitNneu} that the first term on the right vanishes when $\delta \rightarrow 0$ if $\eta <2$. Furthermore, the H\"{o}lder-inequality yields that
\[
	\int_{\Omega \setminus \Omega_+ }\Vert T^+ y - R_{\hat{\alpha}_*} Q Y_{\delta} \Vert^2 d \mathbb{P} 
	\leq \left( \int_{\Omega } \Vert T^+ y - R_{\hat{\alpha}_*} Q Y_{\delta} \Vert^4 d \mathbb{P} \right)^{1/2} \left( \mathbb{P}(\Omega \setminus \Omega_+) \right) ^{1/2}.
\]
Hence, it follows from proposition \ref{prop:Omega-} that for all $\delta < \delta_0$ with $\delta_0 > 0$ sufficiently small it holds with $\eta := 1/s \geq \tfrac{2}{1+2s}$, where $s \in \left(\tfrac{1}{2}, 1\right]$, that
\[
	\sup_{T^+ y \in T_{\varphi}(R)} \int_{\Omega \setminus \Omega_+ }\Vert T^+ y - R_{\hat{\alpha}_*} Q Y_{\delta} \Vert^2 d \mathbb{P} 
	\leq C_{\text{all}} \lceil \delta^{-(2+\eta)} \rceil \exp \left( - \, \tfrac{1}{2} \, C_2 \, \delta^{2-2\eta} \right),
\]
and finally 
\[
	\lim_{\delta \rightarrow 0} \left( \sup_{T^+ y \in T_{\varphi}(R)} \mathbb{E} \left[ \Vert T^+ y - R_{\hat{\alpha}_*} Q Y_{\delta} \Vert^2 \right] \right) = 0,
\]
whicch completes the proof.
\end{proof}

\begin{Rem}[Numerical procedure]\label{rem:Procedure}
The numerical procedure including the estimation of the noise level can be described with the notations of example \ref{ex:konvRVenfWR} as follows:
\begin{center}\begin{tabular}{ll}
\hline
Choose: & $\tau > K_{2,1}; \quad p > 1; \quad q > 1; \quad n \in \mathbb{N}; \quad m \in \mathbb{N}; \quad \epsilon > 0; \quad \hat{\delta}_0 := 0; \quad k:=0;$ \\
Do: & $k := k+1; $ \\
& $ \hat{\delta}_k := \tfrac{1}{2} \, \tau \, n^{-2} \sum_{j=1}^{n} \left( y_{\delta}(j/n) - y_{\delta}((j-1)/n) \right)^2;$ \\
& $\alpha := \hat{\delta}_k^2; $ \\
& $ n := \max \lbrace n(\alpha, \hat{\delta}_k), p \ast n \rbrace ;$\\
While: & $\left( (k<m) \epsilon + (k>m) \max \lbrace \vert \hat{\delta}_k - \hat{\delta}_j \vert, j = k-m,...,k \rbrace > \epsilon \hat{\delta}_k \right) ; \qquad \qquad$\\
Adapt: & $\kappa := \sqrt{m}; \quad n = n(\alpha, \hat{\delta}_k); \quad B := Q_{n} T; \quad x_1 := (\alpha I + B^* B )^{-1} B^* y_{\delta}; \quad k := 0;$\\
Do: & $k := k+1; $ \\
& $ \alpha := q \ast \alpha; $ \\
& $ n := n(\alpha, \hat{\delta}_k); $ \\
& $ B := Q_{n} T;$\\
& $x_k :=  (\alpha I + B^* B )^{-1} B^* y_{\delta};$ \\
While: & $\left( \Vert x_j - x_k \Vert \leq 4 \kappa \delta \sqrt{\psi(\alpha q^{j-k})}, j \leq k \text{ and } \alpha \leq \Vert T \Vert^2 \right) ;$\\
Return: &$x_{k-1};$\\
\hline
\end{tabular}\end{center}
The second part is a modified version of the strategy presented in \cite{MR2240642}.
\end{Rem}


\section{Conclusion}

In this paper we have developed new concepts for the study of statistical inverse problems. The central idea was to link the noise to the asymptotic of the noise level $\delta \rightarrow 0$, varying its probability distribution, which is assumed to be an element of a fixed class $\mathcal{W}$ w.r.t. which the convergence of the considered regularization is required. By means of this approach we were able to disprove the often supposed general transferability of the Bakushinski{\u\i} veto to the stochastical context. 

\vspace{11 pt}
A lot of continuative issues arise out of this result: The estimation of the noise level gained in importance. In particular estimation methods which utilize just one data set are of special interest as the estimate can be incorporated into a regularization method. How does the various parameter choices react to the usage of an estimated noise level and how can we compensate unwanted behaviors? For which other classes of probability distributions does an analog statement to the Bakushinski{\u\i} veto hold and for which ones can we derive counter examples? 


\section*{Acknowledgment}
The author would like to thank Peter Math\'{e}, WIAS Berlin, and Markus Rei\ss, Humboldt-Universit\"{a}t zu Berlin, for helpful discussions.


\end{document}